\documentclass{elsart3-1}

\usepackage{graphicx}

\usepackage{amssymb}
\usepackage{amsmath}

\usepackage[english,francais]{babel}

\newtheorem{e-proposition}[theorem]{Proposition}

\newtheorem{e-definition}[theorem]{Definition\rm}

\newtheorem{theoreme}{Th\'eor\`eme}

\newtheorem{definition}{D\'efinition}

\renewcommand{\theequation}{\arabic{equation}}
\newcommand{\Sc}{\mathcal{S}}
\newcommand{\Dc}{\mathcal{D}}

\newcommand{\Hc}{{\mathcal{H}}}
\newcommand{\Oc}{{\mathcal{O}}}
\newcommand{\Vc}{{\mathcal{V}}}

\newcommand{\Cc}{\mathcal{C}}
\newcommand{\A}{\mathbb{A}}
\newcommand{\C}{\mathbb{C}}
\newcommand{\R}{\mathbb{R}}
\newcommand{\Q}{{\mathbb{Q}}}
\newcommand{\F}{\mathbb{F}}

\newcommand{\Z}{\mathbb{Z}}

\newcommand{\Db}{\mathbb{D}}
\newcommand{\Lb}{\mathbb{L}}

\newcommand{\Str}{\mathrm{Str}}
\newcommand{\Tr}{\mathrm{Tr}}

\newcommand{\Mrm}{\mathrm{M}}
\newcommand{\GL}{\mathrm{GL}}

\def\og{\leavevmode\raise.3ex\hbox{$\scriptscriptstyle\langle\!\langle$~}}
\def\fg{\leavevmode\raise.3ex\hbox{~$\!\scriptscriptstyle\,\rangle\!\rangle$}}

\begin{document}

\begin{frontmatter}




%
\selectlanguage{francais}
\title{Sym\'etries spectrales des fonctions z\^etas}

\vspace{-2.6cm}
\selectlanguage{english}
\title{Spectral symmetries of zeta functions}



\author[authorlabel1]{Fr\'ed\'eric Paugam}
\ead{frederic.paugam@math.jussieu.fr}

\address[authorlabel1]{Institut de math\'ematiques de Jussieu, 175, rue du
Chevaleret, 75013 Paris}

\begin{abstract}

We define, answering a question of Sarnak in his letter to Bombieri
\cite{Sarnak-Bombieri}, a symplectic pairing on the spectral interpretation
(due to Connes and Meyer) of the zeroes of Riemann's zeta function.
This pairing gives a purely spectral formulation of the proof of
the functional equation due to Tate, Weil and Iwasawa, which, in the
case of a curve over a finite field, corresponds to the usual
geometric proof by the use of the Frobenius-equivariant Poincar\'e
duality pairing in etale cohomology. We give another example of a similar construction
in the case of the spectral interpretation of the zeroes of a cuspidal
automorphic $L$-function, but this time of an orthogonal nature. These constructions
are in adequation with Deninger's conjectural program and the arithmetic
theory of random matrices.


\vskip 0.5\baselineskip

\selectlanguage{francais}
\noindent{\bf R\'esum\'e}
\vskip 0.5\baselineskip
\noindent

On d\'efinit, en r\'eponse \`a une question de Sarnak
dans sa lettre a Bombieri \cite{Sarnak-Bombieri},
un accouplement symplectique sur l'interpr\'etation spectrale
(due a Connes et Meyer) des z\'eros de la fonction z\^eta.
Cet accouplement donne une formulation purement spectrale de
la d\'emonstration de l'\'equation fonctionnelle due
\`a Tate, Weil et Iwasawa, qui, dans le cas d'une courbe
sur un corps fini, correspond \`a la d\'emonstration g\'eom\'etrique usuelle
par utilisation de l'accouplement de dualit\'e de Poincar\'e frobenius-\'equivariant
en cohomologie \'etale.
On donne un autre exemple d'accouplement similaire dans
le cas de l'interpr\'etation spectrale des z\'eros de la fonction
$L$ d'une forme automorphe cuspidale, mais cette fois-ci de
nature orthogonale. Ces constructions sont en ad\'equation
avec les pr\'evisions du programme conjectural de Deninger
et de la th\'eorie arithm\'etique des matrices al\'eatoires.


\end{abstract}
\end{frontmatter}


\selectlanguage{francais}

\section{Sym\'etrie spectrale pour la fonction z\^eta}
\label{zeta}
Le lecteur est inform\'e que les notations $\Hc^i$, bien pratiques
pour des raisons d'analogie entre corps de nombres et corps de fonctions,
ne d\'esignent pas ici des espaces de cohomologie.

Soit $K$ un corps global et $\Oc_{K}$ son anneau d'entiers.
Soit $C_{K}=\A_{K}^*/K^*$
sont groupe de classes d'id\`eles et $|.|:C_{K}\to \R_{+}^*$
la norme.
On rappelle maintenant les notations de Meyer dans \cite{Meyer}.
Si $I$ est un intervalle de $\R$, on note $\Sc(C_{K})_{I}$ l'espace des
fonctions $F$ de Bruhat-Schwartz (voir \cite{Osborne}) sur $C_{K}$
dont la transform\'ee de Mellin-Laplace
$$M(F)(s):=\int_{C_{K}}f(x)|x|^sdx$$
converge dans la bande
$\{z|\,\mathrm{Re}(s)\in I\}\subset \C$.
On notera $J(F):=\frac{1}{|x|}F\left(\frac{1}{x}\right)$
pour $F$ une fonction sur $C_{K}$.
Soit $\Sc(C_{K})_{>1<}=\Sc(\R_{+}^*)_{]1,\infty[}\oplus \Sc(\R_{+}^*)_{]-\infty,0[}$.
On plonge $\Sc(C_{K})_{]-\infty,\infty[}$ dans
$\Sc(C_{K})_{>1<}$ par l'application $i:F\mapsto (F,JF)$.
Si $f$ est une fonction de Schwartz sur $\A$, on notera
$$\Sigma(f)(x):=\sum_{q\in K^*}f(qx)$$
la fonction sur $C_{K}$ correspondante.
On d\'efinit ainsi un plongement de $\Sc(\A)$ dans
$\Sc(C_{K})_{>1<}$ par $f\mapsto (\Sigma(f),J\Sigma(f))$.
On note $\Sc(\A)_{0}$ l'espace des fonctions de Schwartz sur
les ad\`eles telles que $\Sigma(f)\in \Sc(C_{K})_{]-\infty,\infty[}$.
L'application
$$\Sigma:\Sc(\A)_{0}\to \Sc(C_{K})_{]-\infty,\infty[}$$
est $\A_{K}^*$-\'equivariante et donne une action
de $C_{K}$ sur le $\C$-espace vectoriel topologique
$$\Hc:=\Sc(C_{K})_{]-\infty,\infty[}/\Sigma(\Sc(\A)_{0}).$$
Si $\A^*_{(1)}\subset \A_{K}^*$ d\'esigne le noyau de la norme,
on notera
$$\Hc^1(\Oc_{K}):=\Hc^{(\A^*_{(1)})}$$
l'espace vectoriel topologique des invariants.
On a une action naturelle de $N:=|\A_{K}^*|$
sur $\Hc^1(\Oc_{K})$. Selon que $K$ est de caract\'eristique
$0$ ou $p$, on a $N=\R_{+}^*$ ou $N=q^\Z$.
L'application $J$ d\'ecrite plus haut donne une application
$\Hc^1(\Oc_{K})\to \Hc^1(\Oc_{K})$ qu'on peut \'ecrire de mani\`ere \'equivariante
par
$$J:\Hc^1(\Oc_{K})^\vee(-1)\to \Hc^1(\Oc_{K}).$$
De plus, la conjugaison complexe $c$ agit sur $\Hc$ (resp. $\Hc^1$) par
$F\mapsto \bar{F}$. 

Le th\'eor\`eme suivant est une variation due \`a Meyer
\cite{Meyer} sur les travaux de Connes \cite{Connes6}.
Sa d\'emonstration se base sur le fait, utilis\'e aussi
par Godement et Jacquet pour des repr\'esentations
automorphes plus g\'en\'erales dans \cite{Godement-Jacquet},
que la fonction $L$ de $K$ peut-\^etre d\'ecrite comme le plus
grand diviseur commun de la famille
des transform\'ees de Mellin $M(f)(s)$ pour
$f\in \Sc(\A)_{0}$.

\begin{theoreme}
L'action de $N$ sur $\Hc^1(\Oc_{K})$ donne une interpr\'etation spectrale
des z\'eros de la fonction z\^eta de Dedekind (compl\'et\'ee) $L_{K}$
de $K$. L'application $J$
identifie spectralement les z\'eros de $L_{K}(|.|^s)$ avec les
z\'eros de $L_{K}(|.|^{1-s})$.
\end{theoreme}

Remarquons que les fonctions $L$ ci-dessus sont de nature diff\'erente
pour les corps de nombres (ce sont des fonctions sur $\C$) et sur les
corps de fonctions (ce sont des fonctions sur $\C/(\frac{2i\pi}{\log q})\Z$).
Afin de rendre la suite ind\'ependante du choix de la caract\'eristique,
nous allons utiliser une construction similaire \`a celle de Deninger
\cite{Deninger2}.
Soit $\Vc$ un espace vectoriel topologique muni d'une action de
$C_{K}$ dont le spectre est discret. On note
$$\Db(\Vc):=\Gamma_{\Cc^\infty}(\Vc\times \R_{+}^*\to \R_{+}^*)^{(C_{K})}$$
l'espace des $C_{K}$-invariants dans les sections $\Cc^\infty$ de la
projection standard (qu'on suppose pour simplifier \`a valeurs nulles en dehors
d'un sous-ensemble fini d'espaces propres),
l'action sur $\R_{+}^*$ \'etant donn\'ee par la norme.
On munit $\Db(\Vc)$ d'une action naturelle de $\R_{+}^*$.
Ceci nous donne un foncteur de ``suspension'' des repr\'esentations de $C_{K}$
vers les modules continus sur $\Lb:=\Db(\C(0))$ (o\`u $\C(0)$ d\'esigne la
repr\'esentation triviale de $N$ sur $\C$) munis d'une
repr\'esentation de $\R_{+}^*$.
On remarque que dans le cas
des corps de nombres, cette op\'eration de suspension
est isomorphe \`a l'identit\'e sur les $\C$-espaces vectoriels de la forme
$\Hc^1(\Oc_{K})$. Elle n'est donc utile que dans le cas des corps de fonctions.

Un op\'erateur $T$ sur un $\Lb$-module topologique est dit tra\c{c}able
si la s\'erie de ses valeurs propres $\lambda\in \Lb$ converge dans
$\Lb$.
Si $c:V\to V$ est une involution d'un $\Lb$-module topologique
et $T:V\to V$ est un op\'erateur tra\c{c}able sur $V$
qui commute \`a $c$. On d\'ecompose
$V=V^{+}\oplus V^{-}\oplus V^{0}$ non canoniquement comme somme
de deux sous-espaces \'echang\'es par l'involution $c$ et d'un
espace laiss\'e stable. La supertrace
de $T$ sur $(V,c)$ est la classe
$$\Str_{c}(T|V):=[\Tr(T|V^{+})-\Tr(T|V^{-})]\in \Lb/c$$
et la trace positive est le scalaire (non canonique)
$$\Tr_{+,c}(T|V):=\Tr(T|V^{+})\in \Lb.$$

On note $\C(1)$ l'espace vectoriel $\C$ muni de l'action de $C_{K}$
par la norme et soit $H^1_{d}(\Oc_{K}):=\Db(\Hc^1(\Oc_{K}))$ et $\Lb(1)=\Db(\C(1))$.
L'op\'erateur de conjugaison complexe $c$ sur $\Hc^1(\Oc_{K})$
induit une involution $c$ sur
$H^1_{d}(\Oc_{K})$ et une d\'ecomposition non canonique
$H^1_{d}(\Oc_{K})=H^1_{d}(\Oc_{K})^{+}\oplus H^1_{d}(\Oc_{K})^{-}\oplus
H^1_{d}(\Oc_{K})^{0}$
(telle que l'espace d'indice positif corresponde aux z\'eros de partie imaginaire
positive et l'espace d'indice nul aux z\'eros r\'eels).
\begin{definition}
Soit $P_{K}(s)$ le polyn\^ome unitaire dont les z\'eros sont les z\'eros
r\'eels de $L_{K}(s)$ avec multiplicit\'es.
On note $L_{K}^{\pm}(s):=\frac{L_{K}(s)}{P_{K}(s)}$
la fonction obtenue \`a partir de $L_{K}(s)$ en lui \^otant ses
z\'eros r\'eels. On notera aussi $H^{1,\pm}_d(\Oc_{K})$ l'espace
spectral correspondant donn\'e par $H^1_{d}(\Oc_{K})/H^1_{d}(\Oc_{K})^{0}$.
\end{definition}

Remarquons que bien que la fonction z\^eta de Riemann n'aie pas de z\'eros r\'eels
(i.e. le polyn\^ome $P_{\Q}(s)$ vaut $1$ et $L_{\Q}^{\pm}(s)=\hat{\zeta}_{\Q}(s)$
est la fonction z\^eta de Riemann compl\'et\'ee), V. Armitage
a montr\'e dans \cite{Armitage} qu'il existe des corps de nombres
dont la fonction z\^eta de Dedekind s'annule en $1/2$ (i.e. tels que
$P_{K}(1/2)=0$).
Le cas des corps de fonctions est \'etudi\'e en d\'etail par Ramachandran dans
\cite{Ramachandran}.

On a une d\'ecomposition \'equivariante
$$H^{1,\pm}_d(\Oc_{K})=H^1_{d}(\Oc_{K})^{+}\oplus H^1_{d}(\Oc_{K})^{-}.$$
On note $C=\left(\begin{smallmatrix}1 & 0\\ 0 & -1\end{smallmatrix}\right).c$ dans l'\'ecriture matricielle pour cette d\'ecomposition.
On utilise les travaux de Meyer \cite{Meyer} qui montrent
que si $F\in \Sc(C_{K})_{]-\infty,\infty[}$, l'op\'erateur
obtenu par convolution de $F$ le long de l'action de
$C_{K}$ sur $\Hc$ est tra\c{c}able. Si $s,t\in H^1_{d}(\Oc_{K})$, on
note
$$
\psi(s,t):=\Str_{c}(s*Jt|H^1_{d}(\Oc_{K}))=
\Tr(s*Jt|H^1_{d}(\Oc_{K})^{+})-\Tr(s*Jt|H^1_{d}(\Oc_{K})^{-}).
$$

Dans le cas des corps de nombres, cette application
s'identifie \`a
$$\psi(F,G)=\Tr(F*JG-G*JF|\Hc^1(\Oc_{K})^{+})$$
pour $F,G\in \Sc(C_{K})_{]-\infty,\infty[}$.

\begin{theoreme}
L'application $\psi$ est bien d\'efinie et induit un accouplement
antisym\'etrique \'equivariant
$\psi:H^{1,\pm}_{d}(\Oc_{K})\times H^{1,\pm}_{d}(\Oc_{K})\to \Lb(1)$,
de surcroit non d\'eg\'en\'er\'e, et qui identifie spectralement les z\'eros de
$L^{\pm}_{K}(s)$ avec ceux de $L^{\pm}_{K}(1-s)$.
De plus, l'accouplement $\psi(.,C.)$ sera sesquilin\'eaire et
d\'efini positif si l'hypoth\`ese de Riemann pour $L^{\pm}_{K}(s)$ est v\'erifi\'ee.
\end{theoreme}

{\bf Preuve:}
La forme $\psi$ passe au quotient par (ce qu'on pourrait noter
abusivement) $\Db(\Sc(\A)_{0})$ car
(voir \cite{Connes-Marcolli-Consani}, Remark 4.18)
la trace est calcul\'ee comme une somme sur les z\'eros des valeurs
des transform\'ees de Mellin, qui sont justement nulles sur les
images des fonctions dans $\Sc(\A)_{0}$. Ceci implique
que l'application est bien d\'efinie.
Le fait que $\psi$ soit antisym\'etrique d\'ecoule de sa d\'efinition
et de l'\'egalit\'e
$$\Tr(s*Jt|H^1_{d}(\Oc_{K})^{-})=\Tr(t*Js|H^1_{d}(\Oc_{K})^{+})$$
qui d\'ecoule de l'\'equation fonctionnelle.
L'\'equivariance d\'ecoule de l'\'equivariance
de $J:\Hc^\vee(-1)\to \Hc$. On peut l'expliquer plus intuitivement
dans le cas d'un corps de nombres en remarquant que si les
z\'eros de $L_{K}$ sont simples et non r\'eels, et si on note $F=s(1)$ et $G=t(1)$,
$$
\psi(s,t)=
\sum_{\rho\in Z_{+}}\left[M(F)(\rho)M(G)(1-\rho)-M(G)(\rho)M(F)(1-\rho)\right]
$$
o\`u $Z_{+}$ sont les z\'eros de $L_{K}(s)$ de partie imaginaire positive
et on voit par changement de variable dans la transform\'ee de Mellin
que l'action de $\R_{+}^*$ sur $M(F)(\rho)M(G)(1-\rho)$ se fait
par $|.|^\rho.|.|^{1-\rho}=|.|$. Le fait que $\psi$ soit non
d\'eg\'en\'er\'ee d\'ecoule du fait que $J$ est un isomorphisme.
L'accouplement $\psi(.,C.)$ est donn\'e par
$$
\psi(s,C.t):=
\Str_{c}(s*J.\left(\begin{smallmatrix}1 & 0\\ 0 & -1\end{smallmatrix}\right).c.t|
H^1_{d}(\Oc_{K}))=
\Tr(s*Jct|H^1_{d}(\Oc_{K})^{+})+\Tr(s*Jct|H^1_{d}(\Oc_{K})^{-}).$$
Si l'hypoth\`ese de Riemann est vraie, les z\'eros v\'erifient
$\overline{1-\rho}=\rho$, donc
$$
\psi(s,Ct)(1)=\sum_{\rho\in Z_{+}}\left[M(F)(\rho)\overline{M(G)}(\rho)+
\overline{M(G)}(1-\rho)M(F)(1-\rho)\right]=\sum_{\rho}M(F)(\rho)\overline{M(G)}(\rho)
$$
et
$$
\psi(s,Cs)(1)=\sum_{\rho}|M(F)(\rho)|^2>0.
$$
On remarque que dans le cas d'un corps de fonctions $K$
sur un corps fini $\F_{q}$ dont la fonction z\^eta n'a pas
de z\'eros r\'eels,
ceci donne une version ind\'ependante de $\ell$ de l'accouplement
de dualit\'e de Poincar\'e frobenius \'equivariant
$$H^{1}_{et}(X,\Q_{\ell})\times H^{1}_{et}(X,\Q_{\ell})\to \Q_{\ell}(1)$$
sur la courbe $X/\F_{q}$ de corps de fonctions $K$,
l'action du frobenius \'etant identifi\'ee par la th\'eorie du corps
de classe r\'esiduel \`a l'action du groupe $N=q^\Z\subset \R_{+}^*$.
La positivit\'e est dans ce cas un r\'esultat d\^u \`a Weil.
Pr\'ecisons que le passage de l'accouplement en cohomologie \'etale
au notre est fourni par un proc\'ed\'e de suspension \`a la Fontaine
d\'ecrit pr\'ecis\'ement par Deninger dans \cite{Deninger2}, section 3.

L'existence de cet accouplement symplectique est une des nombreuses
pr\'evisions de la th\'eorie qui donne le rapport entre les distributions
des valeurs propres des matrices al\'eatoires et les distributions
des z\'eros des fonctions $L$, dont on trouvera un r\'esum\'e dans
\cite{Katz-Sarnak} et dans le survol \cite{Michel}.
Elle est aussi pr\'evue par le programme de Deninger \cite{Deninger1}.

Une construction explicite de cet accouplement a \'et\'e
d\'ecrite, sous l'hypoth\`ese de Riemann g\'en\'eralis\'ee,
par Sarnak dans sa lettre a Bombieri \cite{Sarnak-Bombieri}. Notre formulation
est ind\'ependante de l'hypoth\`ese de Riemann et r\'epond \`a
une question pos\'ee par Sarnak dans ladite lettre.

\section{Sym\'etrie spectrale pour une repr\'esentation automorphe cuspidale}
On donne maintenant un exemple d'interpr\'etation spectrale dans une
situation de dimension sup\'erieure qui vient de la g\'eom\'etrie, en suivant
de pr\`es le travail de Soul\'e \cite{Soule2} et celui de Meyer \cite{Meyer}.
On remarque que nos constructions vont dans le sens du programme de
Deninger et corroborent (sans d\'emonstration)
le fait qu'uniquement trois types de groupes (symplectique, orthogonal et
unitaire) apparaissent dans le rapport entre z\'eros de fonctions $L$ et matrices
al\'eatoires (voir \cite{Michel}):
le groupe symplectique correspond aux motifs absolus d\'ecoup\'es dans des
$H^{2n+1}$ de vari\'et\'es projectives lisses sur $\Q$, le groupe orthogonal
aux motifs absolus d\'ecoup\'es dans des $H^{2n}$ et le groupe unitaire
aux motifs absolus d\'ecoup\'es dans des vari\'et\'es sur des corps imaginaires.

Soit $E$ le mod\`ele de N\'eron sur $\Z$ d'une courbe elliptique sur $\Q$.
Par le r\'esultat fameux de Wiles \cite{Wiles1}
et Breuil, Conrad, Diamond, Taylor \cite{BCDT}, il existe une forme modulaire
cuspidale de poids $2$ telle que
$$L(E,s)=L(f,s),$$
ce qui fait que la fonction $L$ de $E$ est automorphe.

Soit $\pi$ la repr\'esentation automorphe de $\GL_{2}(\A)$
engendr\'ee par $f$. Par construction, elle a un caract\`ere central
trivial. On notera pour chaque $t\in \R_{+}^*$, $G_{t}$ le sous-ensemble
de $\GL_{2}(\A)$ des matrices de d\'eterminant $t$. Soit 
$\Cc_{\pi}$ l'ensemble des coefficients admissibles de $\pi$
pour son action sur $L^2_{0}(\GL_{2}(\Q)\backslash \GL_{2}(\A),1)$
o\`u $1$ d\'esigne le caract\`ere trivial et soit $\Phi$ une fonction de Schwartz
sur $\Mrm_{2}(\A)$. On d\'efinit une fonction sur $\R_{+}^*$ \`a valeurs
complexes par
$$E(\Phi,f)(t)=\int_{G_{t}}\Phi(g)f(g)dg.$$

Si $I$ est un intervalle de $\R$, on note $\Sc(\R_{+}^*)_{I}$ l'espace des
fonctions de Schwartz sur $\R_{+}^*$ dont la transform\'ee de Mellin
converge dans la bande $\{z|\,\mathrm{Re}(s)\in I\}\subset \C$.
Soit $\Sc(\R_{+}^*)_{><}=\Sc(\R_{+}^*)_{]-\infty,0[}\oplus \Sc(\R_{+}^*)_{]2,\infty[}$.
On dispose d'une application $J_{2}:\Sc(\R_{+}^*)_{><}\to \Sc(\R_{+}^*)_{><}$
donn\'ee par $J_{2}(F)(x)=\frac{1}{|x|^2}F(\frac{1}{x})$.
On notera $\Sigma$ l'application
$$\Sigma:\Sc(\Mrm_{2}(\A))\times \Cc_{\pi}\to \Sc(\R_{+}^*)_{><}$$
donn\'ee par
$\Sigma(\Phi,f)=(E(\hat{\Phi},\overset{\vee}{f}),E(\Phi,f))$.
On a aussi une inclusion diagonale de $\Sc(\R_{+}^*)_{]-\infty,\infty[}$ dans
$\Sc(\R_{+}^*)_{><}$ donn\'ee par $F\mapsto (F,F)$.
On notera $\Hc_{-}$ l'image de cette inclusion diagonale et $\Hc_{+}$
l'adh\'erence de l'espace vectoriel engendr\'e par l'image de $\Sigma$.
Soit $\Hc^2_{p}(E)$ le quotient $\Hc_{+}+\Hc_{-}/\Hc_{+}$.
Cette notation est expliqu\'ee par la d\'ecomposition en motifs
absolus (au sens entendu par Manin \cite{Manin4}), la derni\`ere ligne
ayant un sens plus concret d'apr\`es les r\'esultats du premier paragraphe:
$$
\begin{array}{ccccccc}
\Hc^*(E) & = & \Hc^0(E/\Z) & \oplus & \Hc^1(E/\Z) & \oplus & \Hc^2(E/\Z)\\
& = & \Hc^*(\Z) & \oplus & \Hc^2_{p}(E) & \oplus & \Hc^*(\Z)(1)\\
& = & \Hc^0(E)  \oplus \Hc^1(E)\oplus \Hc^2_1(E)& \oplus & \Hc^2_p(E) &\oplus & \Hc^2_3(E)\oplus \Hc^3(E)\oplus \Hc^4(E)\\
& := & \C\oplus \Hc^1(\Z)\oplus \C(1) & \oplus & \Hc^2_{p}(E) & \oplus & \C(1)\oplus \Hc^1(\Z)(1)\oplus \C(2)
\end{array}
$$
qui correspond par d\'efinition \`a la d\'ecomposition en produit de la
fonction $\Lambda(E,s)$ en
$$\Lambda(E,s)=\frac{2\pi}{s}.L_{\Z}(s).\frac{2\pi}{1-s}.\frac{1}{L(E,s)}.
\frac{2\pi}{s-1}.L_{\Z}(s-1).\frac{2\pi}{2-s}.$$
Ainsi, dans ce formalisme de Manin, la fonction $L(E,s)$ correspond \`a
un morceau du $\Hc^2(E)$ not\'e $\Hc^2_{p}(E)$ qu'on obtient en lui
\^otant ses deux $\C(1)$ superflus.
On dispose d'une application surjective $\Sc(\R_{+}^*)_{]-\infty,\infty[}\to \Hc^2_{p}(E)$
et on peut faire agir $\Sc(\R_{+}^*)_{]-\infty,\infty[}$ sur $\Hc^2_{p}(E)$
par convolution. On note
$$\psi:\Sc(\R_{+}^*)_{]-\infty,\infty[}\times \Sc(\R_{+}^*)_{]-\infty,\infty[}\to \C(2)$$
l'application donn\'ee par
$$\psi(F,G):=\Tr_{+,c}(F*J_{2}G+G*J_{2}F|\Hc^2_{p}(E)^{+})$$
o\`u $\Hc^2_{p}(E)^{+}$ d\'esigne la partie correspondant aux z\'eros
de partie imaginaire positive.
Le choix (relativement arbitraire) de cette forme est bas\'e sur l'analogie entre
corps de nombres et corps de fonctions, ainsi que sur l'analogie entre motifs
absolus et structures de Hodge, qui impose (pour des raisons de poids) que
la forme soit sym\'etrique, puisque la fonction $L$ correspondante est associ\'ee
au $H^2$ de la surface arithm\'etique $E$.

\begin{theoreme}
On suppose que $L(E,s)$ ne s'annule pas sur l'axe r\'eel.
L'application $\psi$ induit un accouplement
$\psi:\Hc^2_{p}(E)\times \Hc^2_{p}(E)\to \C(2)$
\'equivariant, sym\'etrique et non d\'eg\'en\'er\'e qui identifie spectralement les z\'eros de
$L(E,s)$ avec ceux de $L(E,2-s)$. Plus g\'en\'eralement
(en combinant ce r\'esultat avec ceux du premier paragraphe),
on obtient une forme bilin\'eaire
$$\phi:\Hc^*(E)\times \Hc^*(E)\to \C(2)$$
non d\'eg\'en\'er\'ee qui identifie spectralement les z\'eros et p\^oles
de $\Lambda(E,s)$ avec ceux de $\Lambda(E,2-s)$. 
\end{theoreme}

On remarque que dans le cas d'une courbe elliptique $E$ sur un corps
de fonctions, on construit \`a la main comme ci-dessus et dans le langage de la
section \ref{zeta} une forme qui est une version ind\'ependante
de $\ell$ de l'accouplement de dualit\'e de Poincar\'e frobenius \'equivariant sur
la cohomologie \'etale de la courbe. On aimerait disposer d'une mani\`ere
plus canonique de mettre en place toutes ces constructions.

\section*{Remerciements}
Je remercie Miguel Bermudez pour d'utiles discussions,
Alain Connes, Katia Consani, Christopher Deninger, Matilde Marcolli et
Ralf Meyer pour des explications sur leurs travaux et en particulier
Alain Connes pour m'avoir indiqu\'e un probl\`eme dans le traitement des corps
de fonctions dans une version pr\'eliminaire de ce travail.
Merci enfin \`a Jean-Fran\c{c}ois Mestre qui m'a fournit des r\'ef\'erences
sur les z\'eros r\'eels.

\bibliographystyle{alpha}
\bibliography{/Users/fpaugam/Documents/travail/fred}

\end{document}